# DIGITAL ERA: Magic Squares and 8<sup>th</sup> May 2010 (08.05.2010)

## INDER JANEJA

Department of Mathematics
Federal University of Santa Catarina
88040-900 Florianópolis, SC, Brazil.
E-mail: ijtaneja@gmail.com
Web-site: www.mtm.ufsc.br/~taneja

## ABSTRACT

*In this short note we have produced different kind of magic squares using digital letter having only the algorisms:* 0, 1, 2, 5 and 8. *The interesting fact in considering these five digits is that the day 8<sup>th</sup> May 2010 also have these ones (08.05.2010). Moreover, the magic squares presented have some interesting properties, such as: they remains the same if we rotate them by 180<sup>o</sup>, or see in the mirror, or see on the other side of the paper, etc. Two palindrome semi-magic squares of order 3x3 are also given. Still, we have considered other dates having four digits.*

## DETAILS

It is well known that there are digits specially used in watches, elevators, etc. These digits are of type:

0, 1, 2, 3, 4, 5, 6, 7, 8 and 9.

There are some of these algorisms digits that have special property such as when we rotate them they become known algorisms again. These are

0, 1, 2, 5, 6, 8 and 9.

In these cases, when we rotate them to 180<sup>o</sup> (degrees) they again remain the known one. This means that the digits

0, 1, 2, 5 and 8,

remains the same while

6 becomes as 9 and 9 as 6.

Summarizing we can say that there only five digits that remains the same when we them gave a 180<sup>o</sup> rotation and these are

0, 1, 2, 5 and 8

Interesting is that the date MAY 8, 2010 have only these five digits, i.e.,

08/05/2010 or 08.05.2010.

Here below we present for the **first time two fantastic magic squares** with the digits:

0, 1, 2, 5 and 8.

The first one is of order 5x5 and has a combination of all the five digits while the second one is of order 4x4 and has the combination of digits:

1, 2, 5 and 8.

See below:

| 00 | 11 | 22 | 88 | 55 |
|----|----|----|----|----|
| 82 | 58 | 05 | 10 | 21 |
| 15 | 20 | 81 | 52 | 08 |
| 51 | 02 | 18 | 25 | 80 |
| 28 | 85 | 50 | 01 | 12 |

and

| 52 | 11 | 85 | 28 |
|----|----|----|----|
| 88 | 25 | 51 | 12 |
| 21 | 82 | 18 | 55 |
| 15 | 58 | 22 | 81 |

The above magic squares have some common interesting properties:
(i) Both of them can be rotated to 180° and still remains magic squares.
(ii) Both of them remains again the magic squares when we change the order of elements, for example, 52 as 25 etc.;
(iii) Both of them remain again magic square, if we see them in the mirror, or reflection in water or on the other side of the glass, or other side of the paper.
(iv) In all the situations the sum of the lines, columns or two main diagonals of both the magic squares always remains 176=2X88=88+88.

Interestingly, 88+88 remains the same under all the above situations. See again the 5x5 order **universal magic square** with border sum:

| 88+88 | 88+88 | 88+88 | 88+88 | 88+88 | 88+88 | 88+88 |
|---|---|---|---|---|---|---|
| 88+88 | 00 | 11 | 22 | 88 | 55 | 88+88 |
| 88+88 | 82 | 58 | 05 | 10 | 21 | 88+88 |
| 88+88 | 15 | 20 | 81 | 52 | 08 | 88+88 |
| 88+88 | 51 | 02 | 18 | 25 | 80 | 88+88 |
| 88+88 | 28 | 85 | 50 | 01 | 12 | 88+88 |
| 88+88 | 88+88 | 88+88 | 88+88 | 88+88 | 88+88 | 88+88 |

*Universal 88+88 magic square*

Until today there is only one magic square of order 4x4 that has all these properties and is only with the digits 1 and 8 (*Ivan Moscovich, Fiendishly Difficult Math Puzzles, Sterling. New York, 1986, p.18*), famous as **IXOHOXI** ou **universal magic square**. Here we presented universal magic square of order 5x5 with five digits combinations and the universal magic square of order 4x4 with four digits combinations. The above 5x5 order magic squares still have the property of pen-diagonal.

Alsol we can have the following 3x3 order **palindromic universal semi-magic squares**:

| 252 | 515 | 121 |
|---|---|---|
| 111 | 222 | 555 |
| 525 | 151 | 212 |

| 282 | 818 | 121 |
|---|---|---|
| 111 | 222 | 888 |
| 828 | 181 | 212 |

The sum of rows and columns of the first one is 888 and the second one is 1221. Here the sum of the diagonals is not the same as of the rows or columns. That's why we call it **Semi-magic Square**. But both of them have the properties given in (i), (ii) and (iii).

Some other dates happing during 2010 having five digits are as follows:

08.05.2010 – 18.05.2010 – 28.05.2010
05.08.2010 – 15.08.2010 – 25.08.2010.

Still in the month of May, 2010, there are days that have less number of digits such as:

01.05.2010 – 02.05.2010 – 05.05.2010
10.05.2010 – 11.05.2010 – 12.05.2010
15.05.2010 – 20.05.2010 – 21.05.2010
22.05.2010 – 25.05.2010.

These days have only four digits:   0, 1, 2 and 5

In the month of January, February, October, November and December, 2010 we can find dates having only these four digits too. Based on these four digits appearing in the above dates we have the following interesting *universal 88 magic square*:

| 88 | 88 | 88 | 88 | 88 | 88 |
|----|----|----|----|----|----|
| 88 | 52 | 11 | 05 | 20 | 88 |
| 88 | 00 | 25 | 51 | 12 | 88 |
| 88 | 21 | 02 | 10 | 55 | 88 |
| 88 | 15 | 50 | 22 | 01 | 88 |
| 88 | 88 | 88 | 88 | 88 | 88 |

*Universal 88 magic square*

For more studies on magic squares, we recommend the readers the following two sites where one can find a good collection of work, papers, books, etc.

1. http://www.multimagie.com/indexengl.htm.
2. http://recmath.org/Magic%20Squares.

-------------------------